\documentclass[final,hidelinks,onefignum,onetabnum]{siamart220329}
\usepackage{graphicx}
\usepackage{amsfonts}
\usepackage{matlab-prettifier}
\usepackage[nocompress]{cite}
\usepackage{upquote}
\usepackage{verbatim}

\title{Computation of 2D Stokes flows via lightning and AAA rational approximation\thanks{Submitted to the editors DATE.
\funding{YX would like to thank financial support from the UK EPSRC (EP/W522582/1). SLW and YX are grateful to funding from the UK MRC (MR/T015489/1).}}}
\author{Yidan Xue\thanks{Mathematical Institute, University of Oxford, Oxford OX2 6GG, UK (\email{xue@maths.ox.ac.uk}, \email{waters@maths.ox.ac.uk}, \email{trefethen@maths.ox.ac.uk}).}
\and Sarah L. Waters\footnotemark[2]
\and Lloyd N. Trefethen\footnotemark[2]
}
\date{March 2023}

\begin{document}

\maketitle

\begin{abstract}
Low Reynolds number fluid flows are governed by the Stokes equations. In two dimensions, Stokes flows can be described by two analytic functions, known as Goursat functions. Brubeck and Trefethen \cite{Brubeck2022} recently introduced a lightning Stokes solver that uses rational functions to approximate the Goursat functions in polygonal domains. In this paper, we present the ``LARS'' algorithm (Lightning-AAA Rational Stokes) for computing 2D Stokes flows in domains with smooth boundaries and multiply-connected domains using lightning and AAA rational approximation \cite{Nakatsukasa2018}. After validating our solver against known analytical solutions, we solve a variety of 2D Stokes flow problems with physical and engineering applications. Using these examples, we show rational approximation can now be used to compute 2D Stokes flows in general domains. The computations take less than a second and give solutions with at least 6-digit accuracy.
\end{abstract} 



\section{Introduction}
At small Reynolds numbers, where viscous forces dominate inertial forces, fluid flows are governed by the Stokes equations. This is typically true when the fluid is highly viscous or the length scale is small. Stokes flows have numerous physical and biological applications, including microcirculation \cite{Secomb2017}, microfluidic devices \cite{Yeo2011}, swimming microorganisms \cite{Lauga2009}, fluid mixing \cite{Finn2001} and lubrication \cite{Ockendon1995}. Many of the flow characteristics in these problems can be analysed in two dimensions.

For two-dimensional Stokes flows, analytical solutions have been derived for certain problems using the Wiener-Hopf method \cite{Wiener1931,Jeong2001}, the method of images \cite{Jeffery1922,Frazer1926,Finn2001}, and expansion of the stream function \cite{Moffatt1964}. Analytical solutions usually only exist for problems with very simple geometries and boundary conditions. For more complex cases, one can approximate the solutions semi-analytically using the unified transform method \cite{Fokas2003,Luca2018a,Luca2018b} and extended lubrication theory \cite{Tavakol2017}, or solve the Stokes equations using numerical methods including finite element methods \cite{Logg2012}, boundary integral methods \cite{Pozrikidis1992}, and the lattice Boltzmann method \cite{Kruger2017}. It should be noted that there is no clear boundary between semi-analytical and numerical methods. For example, it is sometimes necessary to approximate the boundary conditions using orthogonal polynomials when using the unified transform method \cite{Luca2018b}.

Brubeck and Trefethen \cite{Brubeck2022} recently introduced a ``lightning'' Stokes solver, which uses rational functions to approximate the Goursat functions \cite{Goursat1898}, which are two analytic functions that represent Stokes flow in 2D (see \cref{sec:2d_stokes} for details). The lightning solver differs from other numerical methods, since it treats corner singularities of the Stokes equations by clustering the poles of rational functions exponentially nearby. This enables its root-exponential convergence and thus ``lightning'' computation \cite{Gopal2019a}. In \cite{Brubeck2022}, most Stokes flows were computed to at least 8-digit accuracy in less than a second. For the classic lid-driven cavity problem, the lightning solver captured several self-similar Moffatt eddies \cite{Moffatt1964} near the bottom corners, which can be challenging to resolve using more standard discretisation methods. The Goursat representation for 2D Stokes flows and the lightning solver will be reviewed in \cref{sec:2d_stokes,sec:lightning} of this paper respectively.

However, Brubeck and Trefethen \cite{Brubeck2022} primarily focus on 2D Stokes flow problems in polygons. In this paper, we introduce LARS (Lightning-AAA Rational Stokes), a solver that can be implemented easily to compute 2D Stokes flows in general domains with custom boundary conditions in less than a second. LARS uses several rational approximation algorithms: the lightning solver for sharp corners \cite{Gopal2019a,Brubeck2022}, the AAA rational approximation for smooth boundaries \cite{Nakatsukasa2018,Costa2023} and the series method for multiply connected domains \cite{Axler1986,Finn2003,Trefethen2018}.

It is well known that even for regions with analytic boundaries, the Goursat functions may only be analytically continuable a very short distance across the boundary, an effect known as the ``crowding phenomenon'' \cite{Davis1958,Gopal2019b}. This means that poles may need to be placed close to curved boundaries to achieve good rational approximations. Recently Costa and Trefethen \cite{Costa2023} showed that using AAA rational approximation \cite{Nakatsukasa2018} to place poles outside the curved boundary enables fast and effective solution of Laplace problems. The AAA algorithm, derived from ``Adaptive Antoulas-Anderson'', automatically searches for a rational approximation in barycentric form for a vector of boundary values on a given boundary \cite{Nakatsukasa2018}, and has been implemented in the Chebfun toolbox in MATLAB \cite{Driscoll2014}. In this work, we apply AAA rational approximation to compute Stokes flows in domains with curved boundaries. Using an example case of Stokes flows in a channel with a smooth constriction, we compare our solution against a solution approximated using extended lubrication theory \cite{Tavakol2017}, and present these results in \cref{sec:aaa}.

We then introduce an algorithm for computing Stokes flows in multiply connected domains. For Laplace problems in multiply connected domains, the solution can be approximated using a Laurent series with a logarithmic term \cite{Axler1986}. The series method has been applied to compute numerical solution to Laplace problems \cite{Trefethen2018}. It has also been applied to 2D Stokes flows in domains bounded by cylinders \cite{Price2003,Finn2003}. In this paper, we present an algorithm using the series method to compute Stokes flows in general multiply connected domains. We validate the computed stream function for Stokes flows between two cylinders with different boundary conditions against an analytical solution \cite{Finn2001} in \cref{sec:series_method}.

We emphasise that the main contribution of this paper is the development of a new algorithm for solving Stokes flow problems, and a summary of our numerical method is given in \cref{sec:lars}.  In \cref{sec:application_examples}, we apply LARS to compute various 2D Stokes flow problems to demonstrate its broad application. Using these examples, we show that rational approximation can now be used to compute 2D Stokes flows in general domains. The computation usually takes a fraction of a second for a solution to at least 6-digit accuracy.

\section{2D Stokes flow and biharmonic equations}
\label{sec:2d_stokes}
Define $(x,y)$ as the usual Cartesian coordinate system with associated velocity components $(u,v)$. The steady-state Stokes equations in two dimensions are
\begin{align}
    \mu\nabla^2\mathbf{u}=\nabla{p}, \label{eq:stokes1} \\
    \nabla\cdot\mathbf{u}=0, \label{eq:stokes2}
\end{align}
where $\mathbf{u}=(u,v)^T$ is the 2D velocity field, $p$ is the pressure and $\mu$ is the viscosity. We consider 2D Stokes flow problems \cref{eq:stokes1,eq:stokes2} in a bounded domain $\Omega$. Two boundary conditions are imposed on the domain boundary $\partial\Omega$.

Since the flow is 2D and incompressible, a stream function $\psi$ can be defined by
\begin{equation}
    u=\frac{\partial\psi}{\partial{y}},\ v=-\frac{\partial\psi}{\partial{x}}.
\end{equation}
Next we define the vorticity magnitude $\omega$ as
\begin{equation}
    \omega = \frac{\partial{v}}{\partial{x}}-\frac{\partial{u}}{\partial{y}} = -\nabla^2\psi.
\end{equation}
Taking the curl of \cref{eq:stokes1} gives
\begin{equation}
    \nabla^2\omega=0.
\end{equation}
The stream function thus satisfies the biharmonic equation
\begin{equation}
    \nabla^4\psi=0. \label{eq:biharmonic}
\end{equation}

The Stokes problem now becomes that of finding a solution for the biharmonic equation \cref{eq:biharmonic} in the domain of interest, subject to given boundary conditions. In the complex plane $z=x+iy$, where $x,y\in\mathbb{R}$ and $i=\sqrt{-1}$, we have
\begin{equation}
    \frac{\partial}{\partial{z}}=\frac{1}{2}\left(\frac{\partial}{\partial{x}}-i\frac{\partial}{\partial{y}}\right),\ \frac{\partial}{\partial{\bar{z}}}=\frac{1}{2}\left(\frac{\partial}{\partial{x}}+i\frac{\partial}{\partial{y}}\right),
\end{equation}
where $\bar{z}=x-iy$ is the complex conjugate of $z$. \Cref{eq:biharmonic} can then be written in complex form as
\begin{equation}
\frac{\partial^4\psi}{\partial^2z\partial^2\bar{z}}=0,
\end{equation}
which has a solution
\begin{equation}
\psi(z,\bar{z})=\mathrm{Im}[\bar{z}f(z)+g(z)],
\end{equation}
where $f(z)$ and $g(z)$ are two analytic functions, known as Goursat functions \cite{Goursat1898}.

The flow velocity, pressure and vorticity can be expressed in terms of Goursat functions as
\begin{align}
u-iv = -\overline{f(z)}+\bar{z}f'(z)+g'(z), \label{eq:bc1} \\
\frac{p}{\mu}-i\omega=4f'(z), \label{eq:bc2}
\end{align}
where $\overline{f(z)}$ is the complex conjugate of $f(z)$. In the rest of this paper, \cref{eq:bc1,eq:bc2} will be used to determine the Goursat functions by imposing boundary conditions on $\partial\Omega$ (e.g., $-\overline{f(z_0)}+\bar{z_0}f'(z_0)+g'(z_0)=0$ for a zero velocity boundary condition at $z_0$), from which we can then calculate quantities of interest in the simulation domain $\Omega$.

\section{The lightning solver}
\label{sec:lightning}
Our numerical method is developed from the recently introduced ``lightning'' solver for 2D Stokes flow \cite{Brubeck2022}, which uses rational functions with clustered poles near sharp corners to approximate the Goursat functions, and from the related AAA-least squares method for problems with curved boundaries \cite{Costa2023} (see the next section). These methods are based on the general idea of solving the Laplace equation in polygons or curved domains using rational functions \cite{Gopal2019a}. A typical rational function consisting of $m$ poles $\beta_1,...,\beta_m$ and a polynomial of degree $n$ has the form
\begin{equation}
r(z)=\sum_{j=1}^m\frac{a_j}{z-\beta_j}+\sum_{j=0}^n{b_j}z^j,
\label{eq:rational}
\end{equation}
where $a_j$ and $b_j$ are complex coefficients to be determined from the boundary conditions. In \cite{Gopal2019a} the first and second parts of \cref{eq:rational} are called the ``Newman'' and  ``Runge'' terms, respectively. 

For Stokes flow problems, two rational functions, $\hat{f}(z)$ and $\hat{g}(z)$,  are defined for the two Goursat functions, $f(z)$ and $g(z)$:
\begin{align}
\hat{f}(z) &= \sum_{j=1}^m\frac{a_j^f}{z-\beta_j}+\sum_{j=0}^n{b_j^f}z^j,\\
\hat{g}(z) &= \sum_{j=1}^m\frac{a_j^g}{z-\beta_j}+\sum_{j=0}^n{b_j^g}z^j.
\end{align}

Determining these unknown coefficients is a non-linear problem because of the Newman terms of \cref{eq:rational}. However, it becomes a standard linear least-squares problem if we fix the location of poles beforehand. It has been shown in previous work \cite{Gopal2019a,Brubeck2022} that root-exponential convergence can be achieved if the poles are exponentially clustered near each sharp corner of the domain. For a polygonal domain $\Omega$ with $K$ corners $w_1,...,w_K$, we place $N$ poles near each corner using
\begin{equation}
\beta_{kn} = w_k+Le^{i\theta_k}e^{-\sigma(\sqrt{N}-\sqrt{n})},\ k=1,...,K,\ n=1,...,N,
\label{eq:poles}
\end{equation}
where $L$ is the characteristic length scale, $\theta_k$ is the exterior bisector of corner $w_k$ and $\sigma$ is a constant (normally set as 4), as \cite{Gopal2019a,Brubeck2022}. Note that these lightning poles are only used when the domain boundary has sharp corners and they do not appear in smooth boundary problems. As we will show in the next section, the pole vector $\beta$ for smooth boundaries can be obtained easily using the AAA algorithm \cite{Nakatsukasa2018}.

The representation \cref{eq:rational} of a rational function can be ill-conditioned. Here we carry out a Vandermonde with Arnoldi (VA) orthogonalization \cite{Brubeck2021} for the Runge terms and the group of poles near each corner to construct a well-conditioned basis for the linear system. There are two issues to note here. Firstly, unlike Laplace problems, Stokes flow problems involve the derivatives of Goursat functions. These need to be calculated based on the new basis from the VA orthogonalization (see Equations (4.4)--(4.6) in \cite{Brubeck2022}). Secondly, the Laurent series used in multiply connected problems will also need to be orthogonalized, which will be further discussed in Section \ref{sec:series_method}. 

The sample points are selected along the boundary $\partial\Omega$, and are also clustered near the sharp corners \cite{Gopal2019a,Brubeck2022}. Along smooth boundary components, the sample points are evenly distributed, although improvements would certainly be possible here for cases of strong curvature. By applying two boundary conditions at each sample point using \cref{eq:bc1} or \cref{eq:bc2}, we obtain a well-conditioned least-squares problem $Ax\approx b$. The real matrix $A$ has size $2M\times4(m+n+1)$ and the real vector $b$ has size $2M$, where $M$ is the number of sample points. The columns of the matrix $A$ correspond to the real and imaginary parts of the complex coefficients $a_j^f$, $b_j^f$, $a_j^g$ and $b_j^g$, while its rows correspond to the two boundary conditions applied at $M$ sample points, the values of which are stored in the vector $b$. The solution $x$ gives the optimal coefficients for the two rational approximations $\hat{f}(z)$ and $\hat{g}(z)$ (for the two Goursat functions $f(z)$ and $g(z)$), which satisfy the boundary conditions on $\partial\Omega$ in a least-squares sense. The least-squares problem can be solved easily using the backslash command in MATLAB.

\section{AAA rational approximation for a curved boundary}
\label{sec:aaa}

In biological and engineering applications, many boundary components are curved \cite{Lauga2009,Tavakol2017,Secomb2017,Kelley2023}. The solutions associated with such problems are analytic. However, they may be analytically continuable only a very short distance across the boundary, an effect known as the ``crowding phenomenon'' \cite{Davis1958,Menikoff1980,Banjai2008,Gopal2019b}. In such cases, accurate rational approximations will need to have poles very close to the boundary. 

One way to tackle this phenomenon in rational approximation is to use the AAA algorithm \cite{Nakatsukasa2018}. In brief, the AAA algorithm searches for poles for a rational approximation speedily, reliably and automatically. This algorithm has been found to be very fast and effective in rational approximation of conformal maps near singularities \cite{Gopal2019b}. In a recent work, Costa and Trefethen \cite{Costa2023} applied the AAA algorithm to solve Laplace problems. The algorithm was shown to be able to place poles near the boundary of an arbitrary domain in a configuration effective for rational approximation. Using these poles, 8-digit accuracy was easily achieved.

Here we further apply AAA rational approximation to Stokes flow. We demonstrate the outstanding ability of this algorithm to place poles effectively near curved boundaries using an example case of Stokes flow in a channel with a smooth constriction. We choose this case because the pressure drop across the constriction for a given inlet flux has been determined semi-analytically using extended lubrication theory \cite{Tavakol2017}. We first present the problem in \cref{sec:channel_problem}, before computing the problem using polynomials in \cref{sec:polynomial} and rational functions in \cref{sec:aaa_channel}.

\subsection{Stokes flow in a smoothly constricted channel}
\label{sec:channel_problem}

\begin{figure}[htbp]
  \centering
  \includegraphics[width=.7\textwidth]{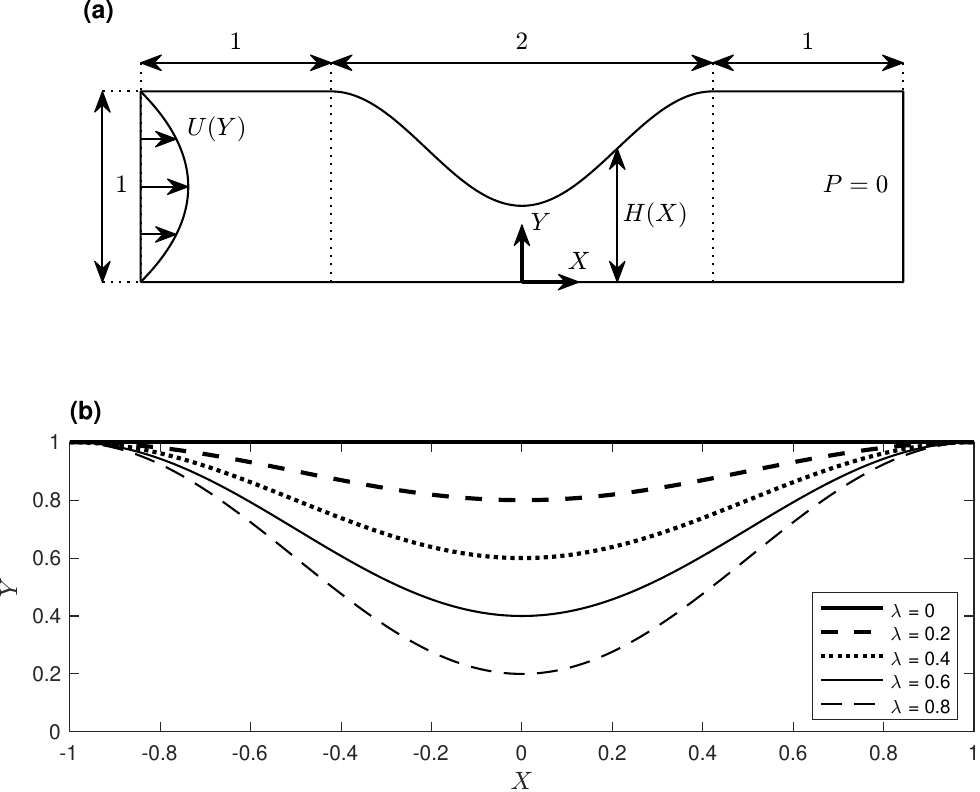}
  \caption{Schematic of Stokes flow through a smoothly constricted channel, after \textnormal{\cite{Tavakol2017}}. (a) Geometry and boundary conditions. (b) Shape function $H(X)$ of the upper boundary for different $\lambda$.}
  \label{fig:constriction}
\end{figure}

We consider Stokes flow through a channel with characteristic length $L_0$ and height $h_0$ with $\delta=h_0/L_0$, and inlet flux $q_0$. We introduce dimensionless variables after \cite{Tavakol2017}:
\begin{equation}
X = \frac{x}{L_0},\ Y = \frac{y}{h_0},\ U=\frac{u}{q_0/h_0},\ V = \frac{v}{q_0/L_0},\ P = \frac{p}{\mu{q_0}L_0/h_0^3},
\end{equation}
where we denote dimensionless variables with capitals and the dimensionless velocity field is $\mathbf{U}=(U,V)^T$. \Cref{eq:stokes1,eq:stokes2} can then be written in the dimensionless form as
\begin{align}
\delta^2\frac{\partial^2U}{\partial{X}^2}+\frac{\partial^2U}{\partial{Y}^2}=\frac{\partial{P}}{\partial{X}}, \label{eq:channel_nd1}\\
\delta^4\frac{\partial^2V}{\partial{X}^2}+\delta^2\frac{\partial^2V}{\partial{Y}^2}=\frac{\partial{P}}{\partial{Y}}, \label{eq:channel_nd2}\\
\frac{\partial{U}}{\partial{X}}+\frac{\partial{V}}{\partial{Y}}=0. \label{eq:channel_nd3}
\end{align}
\Cref{fig:constriction}a presents the Stokes flow problem in a channel with a smooth constriction from $X=-1$ to $X=1$ and $\delta=1$. The shape function for constriction is defined as
\begin{equation}
H(X)=1-\frac{\lambda}{2}(1+\cos(\pi{X})),\ 0\leq\lambda<1,
\label{eq:constriction_shape}
\end{equation}
where $\lambda$ is the maximum dimensionless amplitude of constriction. The boundary curves corresponding to different values of $\lambda$ from 0 to 0.8 are shown in \cref{fig:constriction}b.

A Poiseuille inlet velocity profile $\mathbf{U}(Y)=(6(Y-Y^2),0)^T$, zero velocity on the walls ($\mathbf{U}=\mathbf{0}$), and a parallel outflow profile with zero pressure ($V=0$, $P=0$) are imposed on the domain boundary.

For $0 < \delta \ll 1$, corresponding to a channel with small aspect ratio, we can use lubrication theory \cite{Ockendon1995}. We expand $U$, $V$ and $P$ in terms of $\delta$:
\begin{align}
U(X,Y) &= U_0(X,Y)+\delta^2U_2(X,Y)+\delta^4U_4(X,Y)+\cdots,\\
V(X,Y) &= V_0(X,Y)+\delta^2V_2(X,Y)+\delta^4V_4(X,Y)+\cdots,\\
P(X,Y) &= P_0(X,Y)+\delta^2P_2(X,Y)+\delta^4P_4(X,Y)+\cdots,
\end{align}
and solve \cref{eq:channel_nd1,eq:channel_nd2,eq:channel_nd3} at different orders of $\delta$. The solutions for the pressure drop across the smoothly constricted channel at different orders of $\delta$ have been given in \cite{Tavakol2017} :
\begin{align}
\Delta{P}_0(\lambda) &=\frac{3(3\lambda^2-8\lambda+8)}{(1-\lambda)^{5/2}},\\
\Delta{P}_2(\lambda) &=\frac{12\pi^2\lambda^2}{5(1-\lambda)^{3/2}},\\
\Delta{P}_4(\lambda) &=\frac{8\pi^4(428(-1+\sqrt{1-\lambda})-214(-2+\sqrt{1-\lambda})\lambda-53\lambda^2)}{175\sqrt{1-\lambda}}.
\end{align}
In classical lubrication theory (CLT), only the leading order solution $\Delta{P}_0$ is used for approximating the pressure drop. Tavakol et al. \cite{Tavakol2017} approximate the solution using the leading order term with higher order correction terms, leading to 2nd-order extended lubrication theory (ELT):
\begin{equation}
\Delta{P}(\lambda) =\Delta{P}_0(\lambda)+\delta^2\Delta{P}_2(\lambda)+\mathcal{O}(\delta^4),
\end{equation}
and 4th-order ELT:
\begin{equation}
\Delta{P}(\lambda) =\Delta{P}_0(\lambda)+\delta^2\Delta{P}_2(\lambda)+\delta^4\Delta{P}_4(\lambda)+\mathcal{O}(\delta^6). \label{eq:elt}
\end{equation}
Tavakol et al. \cite{Tavakol2017} show that including higher order terms significantly improves the approximation accuracy of the pressure drop across the constriction for channels with high aspect ratios, e.g. when $\delta$ approaches $1$. In the following sections, we set $\delta=1$ for all our computations, and we compare our results with the 4th-order ELT for the same $\delta$.

\subsection{Polynomial approximation for smooth Stokes flow problems}
\label{sec:polynomial}
We first approximate the Stokes flow problem using the lightning Stokes solver \cite{Brubeck2022}. When used without AAA, the lightning Stokes solver tackles smooth boundary problems by means of the polynomial or ``Runge'' term of \cref{eq:rational}. A recent example is provided as Figure 7.3 of \cite{Brubeck2022}, where 10-digit accuracy is achieved in a smooth bent channel using a polynomial of degree 300.

\begin{figure}[htbp]
  \centering
  \includegraphics[width=0.65\textwidth]{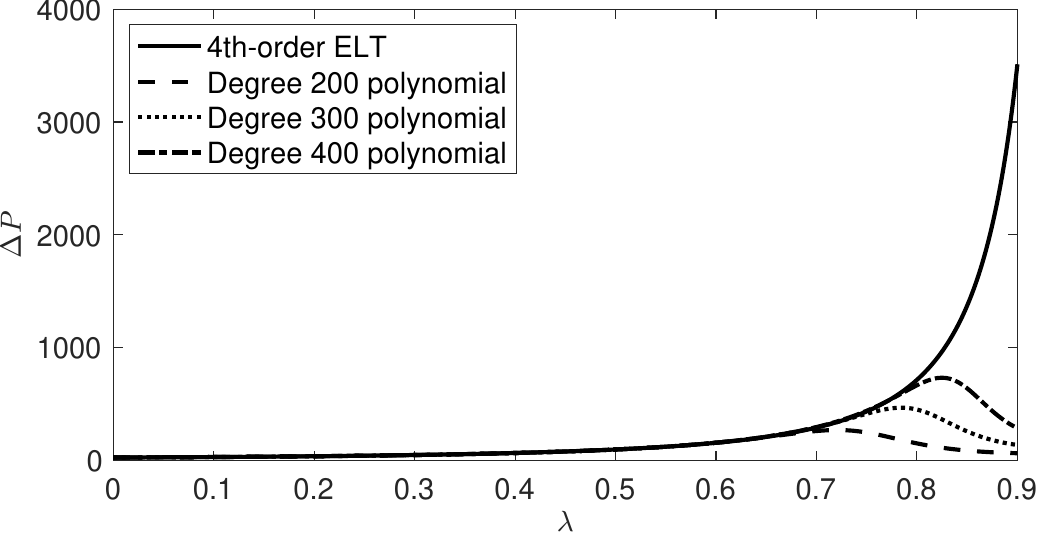}
  \caption{Pressure drop as a function of constriction parameter $\lambda$ when $\delta=1$ computed using polynomials with degrees $200$, $300$ and $400$. The simulation results are compared with the solutions derived using a $4$th-order extended lubrication theory \textnormal{\cite{Tavakol2017}}. The numbers of sample points for the polynomial approximations are $4200$, $6300$ and $8400$.}
  \label{fig:error_polynomial}
\end{figure}

However, the polynomial approximation behaves poorly in this constricted channel problem, especially when the amplitude $\lambda$ is close to 1, due to the crowding phenomenon. \Cref{fig:error_polynomial} presents the pressure drop across the constriction for different $\lambda$ using polynomial approximations of degrees 200, 300 and 400. The numbers of sample points for the approximations are 4200, 6300 and 8400, scaled with the polynomial degree. For example, when the polynomial degree is 200, there are 600 points evenly distributed on each segment of the domain boundary. We treat the upper boundary as 4 boundary segments $X\in(-2,-1)$, $X\in(-1,0)$, $X\in(0,1)$ and $X\in(1,2)$, where 600 points are sampled on each segment. In the complex plane $Z=X+iY$, the pressure drop is calculated between $Z=-1+0.5i$ and $Z=1+0.5i$ using \cref{eq:bc2}.

In \cref{fig:error_polynomial}, the polynomial approximation is compared with the semi-analytical solution approximated by 4th-order ELT \cite{Tavakol2017}. For $\lambda\leq0.7$, the polynomial approximation provides a reasonable estimate of the pressure drop across the constriction (compared with the 4th-order analytical approximation). However, the degree 200 polynomial fails to approximate the pressure drop for $\lambda>0.7$, followed by the degree 300 polynomial for $\lambda>0.75$ and the degree 400 one for $\lambda>0.8$. Note that the computational cost rises sharply as we increase the polynomial degree from 200 to 400. Further increasing the degree of the polynomials may approximate our problem at a larger $\lambda$, but this is certainly not a practical method for all $0\leq\lambda<1$. This is where the AAA algorithm comes into play.

\subsection{AAA rational approximation for the upper boundary}
\label{sec:aaa_channel}
The key idea now is to use the AAA algorithm to place poles outside the curved boundary (which is the upper boundary of the channel in this case) \cite{Nakatsukasa2018} to help the lightning solver compute the Stokes flow \cite{Brubeck2022}. This method has proven to be very effective for solving Laplace problems \cite{Costa2023}. The AAA rational approximation can be computed using the MATLAB code \lstinline[style=Matlab-bw]{aaa.m} in Chebfun \cite{Driscoll2014}, and the workflow is simple:
\vspace{6pt}
\begin{enumerate}
    \item Create a vector $Z_b$ of sample points along the curved boundary.
    \item Apply a boundary function $F$ to $Z_b$, e.g., the Schwarz function \cite{Davis1958}: $F=\overline{Z_b}$.
    \item Run the AAA algorithm for these sample points and boundary values.
    \item Remove the poles inside the domain $\Omega$.
\end{enumerate}
\vspace{6pt}
In MATLAB, the poles can be obtained easily by executing
\vspace{6pt}

\begin{small}
\begin{verbatim}
    F = conj(Zb);
    [r,pol] = aaa(F,Zb,'tol',1e-8);
    inpoly = @(z,w) inpolygon(real(z),imag(z),real(w),imag(w));
    jj = inpoly(pol,Z);
    beta = pol(~jj);
\end{verbatim}
\end{small}
\vspace{6pt}
where $Z_b$ is the vector of sample points along the curved boundary we aim to place poles near, and $Z$ is the vector of sample points along all of $\partial\Omega$. For AAA-lightning computations of the constricted channel problem, we sample 600 points clustered towards singularities at two ends using \lstinline[style=Matlab-bw]{tanh(linspace(-14,14,600))} on each boundary segment.

Here we use the Schwarz function \cite{Davis1958} as the boundary function and set the tolerance of the AAA algorithm as $10^{-8}$ for fast computation. Our numerical experiments show that the use of different boundary functions has negligible effect on the location of poles for this geometry. Here we use the Schwarz function, because it is purely based on the boundary geometry instead of the boundary data, so the poles of the rational approximation of the Schwarz function capture the singularities of the boundary geometry. The pole vector \lstinline[style=Matlab-bw]{beta} will be used to construct the Newman part of our rational function \cref{eq:rational}. For better numerical stability, VA orthogonalization \cite{Brubeck2021} was performed on these poles and the polynomial using \lstinline[style=Matlab-bw]{VAorthog(Z,n,beta)}. This is the same as in the original lightning Stokes solver \cite{Brubeck2022} except that the vector of poles is obtained using AAA. We will call these poles placed by the AAA algorithm ``AAA poles'' in the rest of this paper.

\begin{figure}[htbp]
  \centering
  \includegraphics[width=\textwidth]{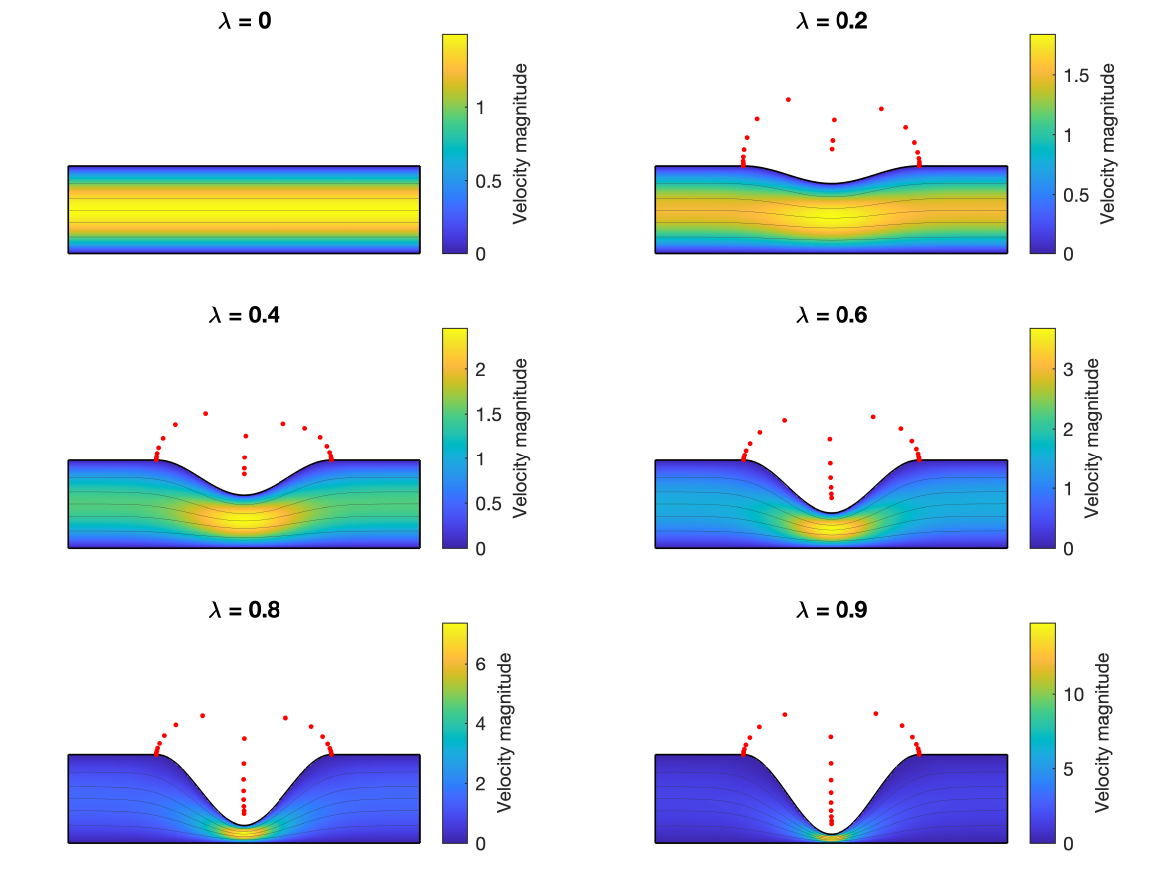}
  \caption{Stokes flow in a smoothly constricted channel for different $\lambda$ from $0$ to $0.9$ and $\delta=1$. The solution is computed using the lighting solver with a polynomial of degree $100$ with poles placed by the AAA algorithm. The locations of poles are marked by red dots. The streamlines and velocity magnitude in each case are represented by light grey lines and a colour scale, respectively. Note that the colour scale has a different range for each case, scaled by the maximum velocity magnitude. The computation for each case takes $1$--$2$ seconds on a standard laptop.}
  \label{fig:channel_cases}
\end{figure}

\Cref{fig:channel_cases} presents the Stokes flow in the constricted channel for different $\lambda$ computed using the algorithm described above. The streamlines in each case are represented by light grey lines, while a colour scale shows the velocity magnitude. The AAA poles are marked as red dots. Using these poles and a polynomial of degree 100, a 5 to 6-digit accurate solution can be computed in 1-2 seconds on a standard laptop. For each computation, running the AAA algorithm, constructing the rational function basis, and solving the least-squares problem takes about 0.8, 0.4 and 0.2 second, respectively. Note that only 3-digit accuracy is achieved in the last case, with $\lambda=0.9$, although the streamlines are qualitatively promising. For higher accuracy, higher degrees of polynomials and more sample points are required. In this paper, ``$\alpha$-digit accuracy'' means that the maximum error of the approximation on the domain boundary is below $10^{-\alpha}$.

In \cref{fig:channel_cases}, the AAA algorithm places poles vertically along the centreline of the constriction, with clustering near the bottom. This phenomenon is very similar to that of the poles placed near a sharp corner using \cref{eq:poles} in the original lightning Stokes solver \cite{Brubeck2022}. However, for this specific problem, it is interesting to note that the AAA poles are clustered towards a point slightly above the boundary instead of on the boundary. This presumably corresponds to a branch point of the analytic continuation across the boundary.

In addition, AAA places poles near the juncture points where the horizontal part transitions into the curved part. The AAA algorithm detects these points as singularities, because the transitions are not twice continuously differentiable (although the first derivative is continuous, see \cref{eq:constriction_shape}).

\begin{figure}[htbp]
  \centering
  \includegraphics[width=0.65\textwidth]{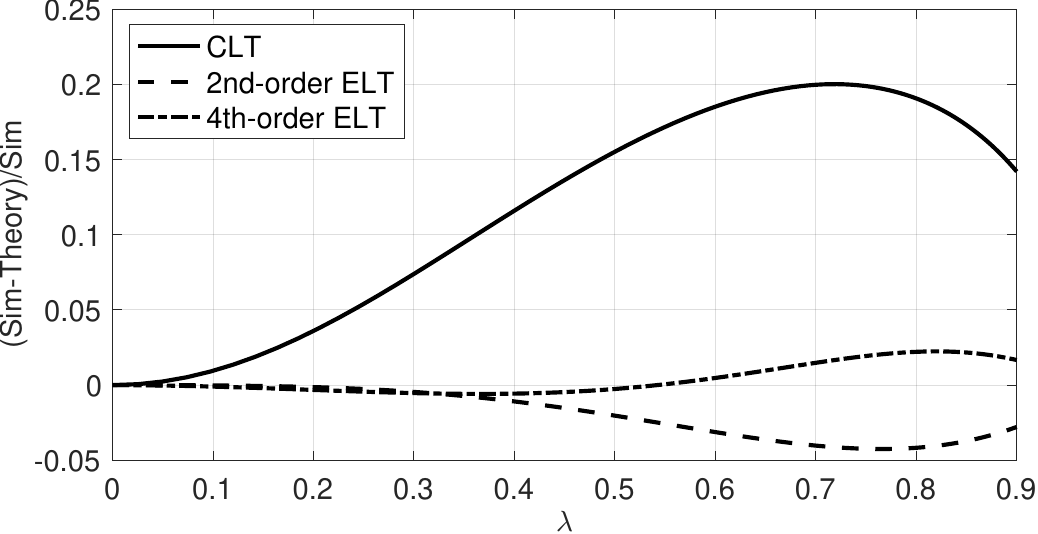}
  \caption{Relative differences between the pressure drop across the constriction when $\delta=1$ approximated by lubrication theory at different orders \textnormal{\cite{Tavakol2017}} and the AAA-lightning computation, which we presume is effectively exact for the purposes of this comparison. The computation uses a degree $100$ polynomial with poles placed by the AAA algorithm.}
  \label{fig:aaa_comparison}
\end{figure}

\Cref{fig:aaa_comparison} shows a comparison of pressure drop across the constriction approximated by lubrication theory at different orders of $\delta$ and the AAA-lightning computation. The gap between lubrication theory and lightning simulation becomes smaller when higher order terms are included in the asymptotic analysis. Compared with our computation, the maximum differences in $\Delta{P}$ for all $\lambda$ are approximately 20\% and 4\%, respectively, when using CLT and 2nd-order ELT. These figures agree with the gaps between the approximations and computations using a finite element method reported in \cite{Tavakol2017}. For 4th-order ELT, the maximum difference in $\Delta{P}$ reduces to 2.2\%.

\section{Series methods for multiply connected domains}
\label{sec:series_method}

In the previous paper \cite{Brubeck2022}, application of the lightning Stokes solver was focused on simply connected domains. However, there has been increasing interest in 2D Stokes flows in multiply connected domains \cite{Ranger1961,Dvinsky1987a,Dvinsky1987b,Manga1996,Finn2001,Price2003,Finn2003,Luca2018a}. In this section, we describe a method to extend the lightning Stokes solver to multiply connected problems.

When solving a Laplace problem in a multiply connected domain, it is known that only a logarithmic term and a Laurent series are needed for each smooth hole \cite{Trefethen2018}. The logarithmic term prevents the solution from being multi-valued, thanks to the ``logarithmic conjugation theorem'' presented by Axler \cite{Axler1986}, with its proof dating to Walsh in 1929 \cite{Walsh1929}. The series method has been shown to be effective for computing conformal maps in multiply connected domains \cite{Trefethen2020}. 

Similar series methods have also been applied to 2D Stokes flow problems in multiply connected domains by Price et al. \cite{Price2003} and Finn et al. \cite{Finn2003}. However, these applications are limited to 2D Stokes flows in domains bounded by cylinders. With the AAA algorithm \cite{Nakatsukasa2018}, VA orthogonalization \cite{Brubeck2021} and the lightning solver \cite{Gopal2019a,Brubeck2022}, we are now able to compute 2D Stokes flows in more general multiply connected domains and thus address much broader applications.

\subsection{Algorithm}

For multiply connected domain problems with $p$ smooth holes, we define a rational function (before adding the logarithmic terms) in the form
\begin{equation}
r(z)=\sum_{j=1}^m\frac{a_j}{z-\beta_j}+\sum_{j=0}^n{b_j}z^j+\sum_{i=1}^p\sum_{j=1}^q{c_{ij}}(z-z_i)^{-j},
\label{eq:rational_multi1}
\end{equation}
where the first and second parts are the Newman and Runge terms from \cref{eq:rational}. The third part represents a Laurent series expansion to degree $q$ at point $z_i$ in the $i$th hole. 

Based on the logarithmic conjugation theorem \cite{Axler1986}, a logarithmic term is needed for each Goursat function. Moreover, unlike Laplace problems, a derivative term $g'(z)$ also appears in \cref{eq:bc1}, which describes the velocity field in Stokes flow problems. The velocity components $u$ and $v$ are each expressed as
\begin{align}
u(z) = \mathrm{Re}[-f(z)+\bar{z}f'(z)+g'(z)], \label{eq:u} \\
v(z) = \mathrm{Im}[-f(z)-\bar{z}f'(z)-g'(z)]. \label{eq:v}
\end{align}
To ensure the velocity field is not multi-valued, an extra term $(z-z_i)\log(z-z_i)-z$ is required in $g(z)$ for the $i$th hole. So the imaginary part of $\log(z-z_i)$ term in $g'(z)$ and $f(z)$ can cancel out in \cref{eq:u,eq:v} \cite{Price2003}. The rational representation of Goursat functions can thus be written as
\begin{align}
\hat{f}(z)&=\sum_{j=1}^m\frac{a_j^f}{z-\beta_j}+\sum_{j=0}^n{b_j^f}z^j+\sum_{i=1}^p\sum_{j=1}^q{c_{ij}^f}(z-z_i)^{-j}+\sum_{i=1}^p{d_i^f}\log(z-z_i),
\label{eq:rational_fz}\\
\hat{g}(z)&=\sum_{j=1}^m\frac{a_j^g}{z-\beta_j}+\sum_{j=0}^n{b_j^g}z^j+\sum_{i=1}^p\sum_{j=1}^q{c_{ij}^g}(z-z_i)^{-j}+\sum_{i=1}^p{d_i^g}\log(z-z_i)\notag\\
&-\sum_{i=1}^p\overline{d_i^f}[(z-z_i)\log(z-z_i)-z],
\label{eq:rational_gz}
\end{align}
where $a_j^f$, $b_j^f$, $c_{ij}^f$, $d_i^f$, $a_j^g$, $b_j^g$, $c_{ij}^g$ and $d_i^g$ are complex coefficients to be determined by solving a least-squares problem.

In the lightning Stokes solver, VA orthogonalization has been carried out for both the Newman and the Runge parts. The VA orthogonalization for the polynomial part can be found in the Section 4.1 of \cite{Brubeck2022}. The VA orthogonalization for the Laurent series is very similar to that for polynomials, and can be realised by appending a new module at the end of the existing MATLAB code \lstinline[style=Matlab-bw]{VAorthog} (see Appendix A in \cite{Brubeck2022}):
\vspace{6pt}

\begin{small}
\begin{verbatim}
    Q = ones(M,1); H = zeros(nl+1,nl);
    for k = 1:nl
       q = 1./(Z-Z_i).*Q(:,k);
       for j = 1:k, H(j,k) = Q(:,j)'*q/M; q = q - H(j,k)*Q(:,j); end
       H(k+1,k) = norm(q)/sqrt(M); Q(:,k+1) = q/H(k+1,k);
    end
    Hes{length(Hes)+1} = H; R = [R Q(:,2:end)];
\end{verbatim}
\end{small}
\vspace{6pt}
where \lstinline[style=Matlab-bw]{nl} is the degree of Laurent series (i.e. $q$ in \cref{eq:rational_fz,eq:rational_gz}), \lstinline[style=Matlab-bw]{Hes} is the upper-Hessenberg matrix from the Arnoldi process. Similarly, the new VA basis can be constructed using \lstinline[style=Matlab-bw]{VAeval} with an additional module appended in the end:
\vspace{6pt}

\begin{small}
\begin{verbatim}
    H = Hes{1}; Hes(1) = []; 
    Q = ones(M,1); D = zeros(M,1);
    Zpki = 1./(Z-Z_i); Zpkid = -1./(Z-Z_i).^2;
    for k = 1:nl
        hkk = H(k+1,k);
        Q(:,k+1) = (Q(:,k).*Zpki - Q(:,1:k)*H(1:k,k))/hkk;
        D(:,k+1) = (D(:,k).*Zpki - D(:,1:k)*H(1:k,k) + Q(:,k).*Zpkid)/hkk;
    end
    R0 = [R0 Q(:,2:end)]; R1 = [R1 D(:,2:end)];
\end{verbatim}
\end{small}
\vspace{6pt}
where \lstinline[style=Matlab-bw]{R0} is the function basis for the Goursat functions $f(z)$ and $g(z)$ and \lstinline[style=Matlab-bw]{R1} is the function basis for their derivatives $f'(z)$ and $g'(z)$. The first column representing the constant term is always omitted at the storing step, since the constant term has been included in the polynomial part. The logarithmic terms are added in the MATLAB code \lstinline[style=Matlab-bw]{makerows}, when evaluating the function values on the domain boundary $\partial\Omega$. Example codes for computing Stokes flow in general multiply connected domains can be found in \cref{sec:code}.

\subsection{Stokes flow between two translating and rotating cylinders}
\label{sec:cylinders}

\begin{figure}[htbp]
  \centering
  \includegraphics[width=0.5\textwidth]{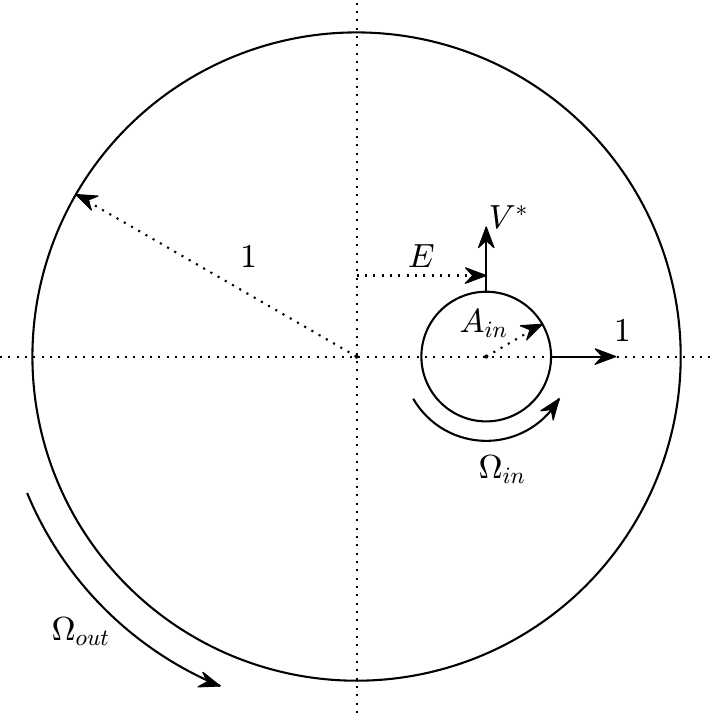}
  \caption{Schematic of a translating and rotating cylinder in a rotating cylinder, after \textnormal{\cite{Finn2001}} (dimensionless quantities).}
  \label{fig:cylinder}
\end{figure}

Here we consider Stokes flow between two cylinders to illustrate the speed and accuracy of our new Stokes flow solver via comparison of numerical results with Finn and Cox \cite{Finn2001}, based on previous work on cylinders with simpler motions \cite{Frazer1926,Ballal1976,Ranger1980}. \Cref{fig:cylinder} shows the problem setting, where the outer cylinder has a radius of $a_{out}$ centred at the origin and the inner cylinder has a radius of $a_{in}$ centred at $(\epsilon,0)$. The outer and inner cylinders rotate with angular velocities $\omega_{out}$ and $\omega_{in}$, respectively. The inner cylinder can also translate with velocity $(u^*,v^*)^T$. We describe the problem using dimensionless variables:
\begin{align}
 (u,v)=u^* (U,V),\ \  p=\frac{u^* a_{out}}{\mu}P ,
\end{align}
where capitals denote dimensionless variables. The system is then characterised by the following dimensionless parameters

\begin{align} V^* = \frac{v^*}{u^*},\ A_{in} = \frac{a_{in}}{a_{out}},\ E = \frac{\epsilon}{a_{out}},\ (\Omega_{in}, \Omega_{out}) = \frac{a_{out}}{u^*} (\omega_{in}, \omega_{out}).
\end{align}
The analytical solution to this problem can be derived by superposing solutions corresponding to each type of motion (rotational and translational), see Equation (68) in \cite{Finn2001} (the details are omitted here). 

\begin{figure}[htbp]
  \centering
  \includegraphics[width=.8\textwidth]{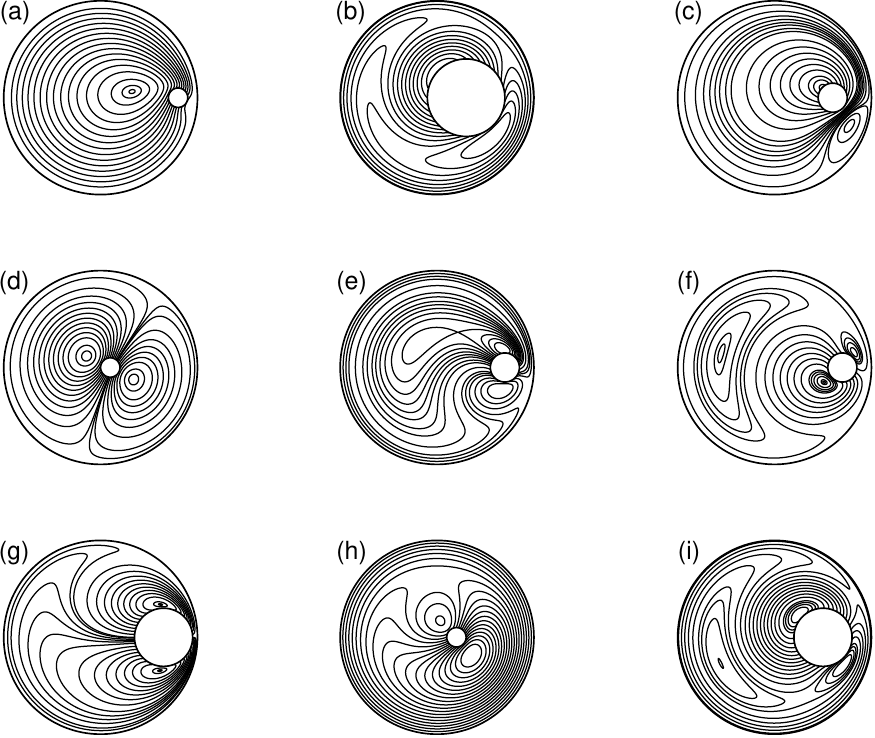}
  \caption{Streamlines for Stokes flow between two cylinders for nine different boundary conditions, following \textnormal{\cite{Finn2001}}. The parameter values are listed in \cref{tab:parameter_value}. The stream function is $0$ on the outer cylinder.}
  \label{fig:finn_cox_sim}
\end{figure}

We used our algorithm to compute the Stokes flow in all nine example cases presented in Figure 9 of \cite{Finn2001}. Our results are shown in \cref{fig:finn_cox_sim}. The parameter values for each case are listed in \cref{tab:parameter_value}. The same parameter values were used in \cite{Finn2001}, except that we change $E=0.7$ to $E=0.65$ in case `g' to prevent the boundaries of the two cylinders touching each other. 

\begin{table}[tbhp]
\footnotesize
  \caption{Parameter values for nine example cases as presented in \cref{fig:finn_cox_sim}.}\label{tab:parameter_value}
\begin{center}
  \begin{tabular}{cccccc} \hline
   Case & $A_{in}$ & $E$ & $V^*$ & $\Omega_{in}$ & $\Omega_{out}$ \\ \hline
    a & 0.1 & 0.8 & 2 & $-3$ & 1 \\
    b & 0.4 & 0.3 & 1 & 5 & $-3$ \\
    c & 0.15 & 0.6 & 1 & 10 & 0 \\
    d & 0.1 & 0.1 & 2 & 0 & 0 \\
    e & 0.15 & 0.7 & 0 & 3.33 & 0.66 \\
    f & 0.15 & 0.7 & $-1$ & 0 & 0.66 \\
    g & 0.3 & 0.65 & 0 & 0 & $-0.2$ \\
    h & 0.1 & 0.2 & 1 & 5 & $-1$ \\
    i & 0.3 & 0.5 & 1 & 2 & $-2$\\ \hline
  \end{tabular}
\end{center}
\end{table}

For each case, 100 sample points are evenly distributed along the inner cylinder boundary, with another 500 points along the outer cylinder boundary slightly clustered towards the narrower gap between two cylinders:
\vspace{6pt}

\begin{small}
\begin{verbatim}
    pw = ceil(1/(1-epsilon))+1;
    sp = tanh(linspace(-pw,pw,500));
    Z_b = a_out*exp(1i*pi*(sp-1)');
\end{verbatim}
\end{small}
\vspace{6pt}
The solution is computed using a rational approximation \cref{eq:rational_multi1} consisting of a degree 20 polynomial and two degree 50 Laurent series about $(E,0)$ and $(1/E,0)$. The computation for each case takes tens of milliseconds on a standard laptop. 

Here we use an additional Laurent series about the reflection of the centre of inner cylinder to better compute the flow field in the narrow gap between two cylinders in cases `a' and `g'. The algorithm using two Laurent series is found to be much more efficient than using a Laurent series about the inner cylinder with a much higher degree polynomial. When the inner cylinder is very close to the outer cylinder, the analytic continuation of Goursat functions may have branch points outside the outer cylinder, leading to the crowding phenomenon that we have tackled using the AAA algorithm in the last section. Hence one may develop another algorithm that places AAA poles near the gap outside the domain. However, controlling the region of AAA poles and dealing with the logarithmic terms can be challenging and requires future work.

\begin{figure}[htbp]
  \centering
  \includegraphics[trim={1.2cm 1cm 1.2cm 0},clip,width=\textwidth]{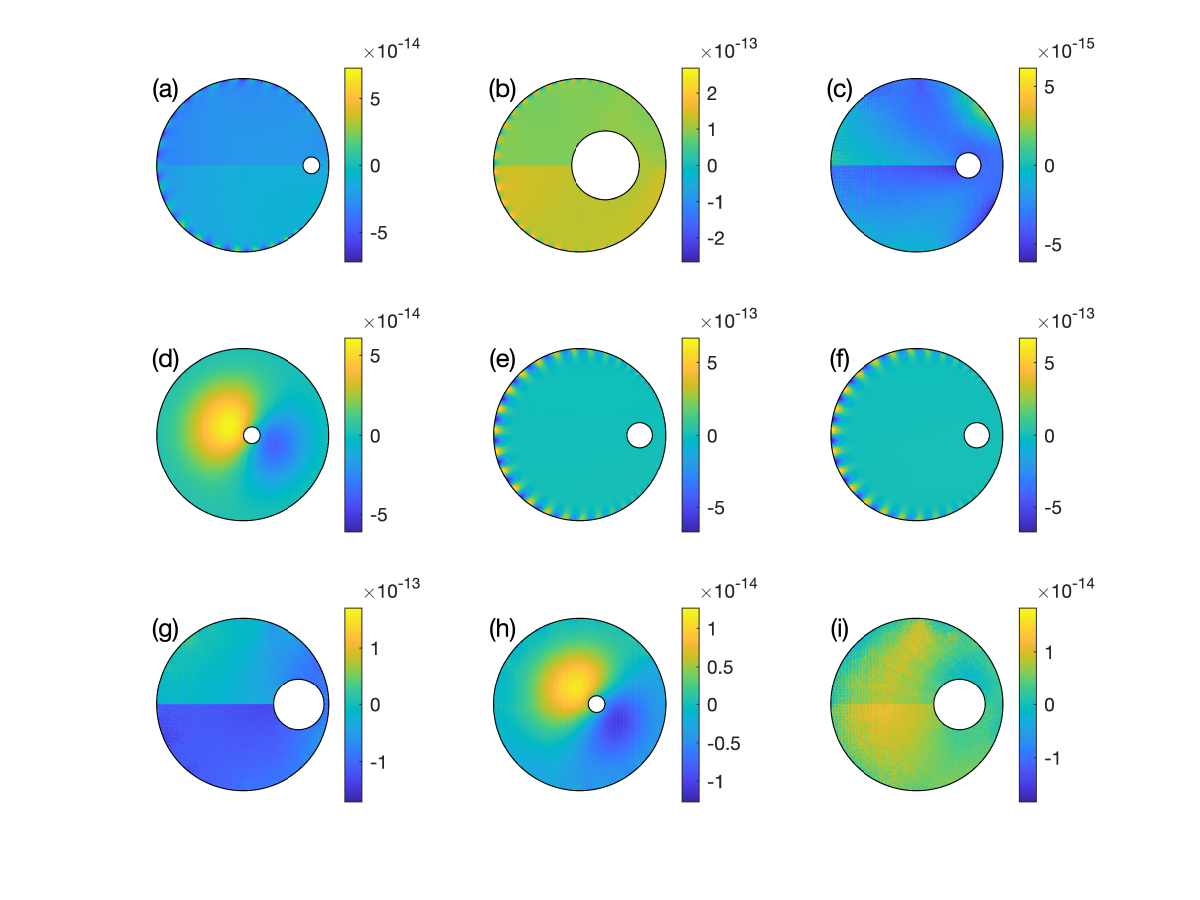}
  \caption{Pointwise error of the stream function computed using our Stokes flow solver in all nine cases.}
  \label{fig:finn_cox_val}
\end{figure}

Next, we validate the values of the stream function $\psi$ computed using our Stokes flow solver against the analytical solution \cite{Finn2001}. \Cref{fig:finn_cox_val} shows the pointwise error ($\mathrm{Error}=\mathrm{Sim}-\mathrm{Theory}$) of our computation of the stream function in all nine cases. All solutions are obtained to 12 to 14-digit accuracy, which is close to machine precision.

\section{Algorithm summary and applications}
\label{sec:solver_application}

In this section, we first summarise the LARS algorithm, and then apply it to a variety of 2D Stokes flow problems to demonstrate the utility of the algorithm.

\subsection{The LARS algorithm} 
\label{sec:lars}

There are three steps to compute the solution to a 2D Stokes flow problem using the LARS algorithm:

\medskip{\em 1. Place the poles and sample points.} For a given domain $\Omega$, we first identify corner singularities and cluster lightning poles exponentially nearby following the lightning algorithm \cite{Gopal2019a,Brubeck2022}. We then identify smooth boundaries and approximate the Schwarz function on each section using the AAA algorithm \cite{Nakatsukasa2018}. After removing the poles inside $\Omega$, we obtain the AAA poles for approximating these smooth boundaries \cite{Costa2023}. The lightning poles and the AAA poles form the vector $\beta_j$ to construct the Newman part in \cref{eq:rational_fz,eq:rational_gz}. The sample points are clustered towards corner singularities along the domain boundary $\partial\Omega$. For smooth domain boundaries and holes, we take equispaced sample points. All sample points along the entire $\partial\Omega$ are stored in a vector $Z=(Z_1,Z_2,...,Z_M)^T$.

\medskip{\em 2. Construct the rational function bases.} We perform VA orthogonalization \cite{Brubeck2021} for the Newman terms, the Runge terms, and the Laurent series corresponding to each hole for all sample points $Z$. Using the VA orthogonalization we obtain a well-conditioned rational function basis $R_0$ (spanning the same spaces as the original basis) to evaluate every sample point in $Z$, and its derivative $R_1$. The matrices $R_0$ and $R_1$ have size $M\times(m+n+pq+1)$, where each row corresponds to a sample point and each column corresponds to a coefficient $a_j$, $b_j$ or $c_{ij}$.

\medskip{\em 3. Solve the least-squares problem and compute physical quantities.} We now impose two boundary conditions on $\partial\Omega$ (with $Z$ being the vector of sample points) to compute the coefficient values in two rational functions $\hat{f}(z)$ and $\hat{g}(z)$ that approximate the Goursat functions. For simplicity, it is assumed here that we have the boundary condition $(\tilde{u}(Z_i),\tilde{v}(Z_i))^T$ at each sample point $Z_i$ on $\partial\Omega$. One can impose $\psi$, $p$, $\omega$ or a component in the 2D stress tensor easily after minor changes from the linear system presented below. Based on \cref{eq:u,eq:v}, we construct a linear system $Ax\approx b$ using $R_0$, $R_1$ and logarithmic terms, which will be solved using a least-squares approach to find the coefficient vector $x$ that minimises $||Ax-b||_2$. The matrix $A$ consists of $2\times8$ blocks:
\begin{align}
A=&\left[
\begin{matrix}
\mathrm{Re}\{cZ{\times}R_1-R_0\} & \mathrm{Re}\{R_1\} & \mathrm{Re}\{cZ{\times}oZ-2lZ\} & \mathrm{Re}\{oZ\}
\\
\mathrm{Im}\{-cZ{\times}R_1-R_0\} & \mathrm{Im}\{-R_1\} & \mathrm{Im}\{-cZ{\times}oZ\} & \mathrm{Im}\{-oZ\}
\end{matrix}
\right.\nonumber
\\
&\qquad\left.
\begin{matrix}
-\mathrm{Im}\{cZ{\times}R_1-R_0\} & -\mathrm{Im}\{R_1\} & -\mathrm{Im}\{cZ{\times}oZ\} & -\mathrm{Im}\{oZ\}
\\
\mathrm{Re}\{-cZ{\times}R_1-R_0\} & \mathrm{Re}\{-R_1\} & \mathrm{Re}\{-cZ{\times}oZ-2lZ\} & \mathrm{Re}\{-oZ\}
\end{matrix}
\right],
\label{eq:lars_linear}
\end{align}
where we have
\begin{align}
cZ = \begin{bmatrix}
\overline{Z_1} & & & \\
& \overline{Z_2} & & \\
& & \ddots & \\
& & & \overline{Z_M}
    \end{bmatrix},\ 
oZ = \begin{bmatrix}
\dfrac{1}{Z_1-z_1} & \dfrac{1}{Z_1-z_2} & \dots & \dfrac{1}{Z_1-z_p} \\
\dfrac{1}{Z_2-z_1} & \dfrac{1}{Z_2-z_2} & \dots & \dfrac{1}{Z_2-z_p} \\
\vdots & \vdots & \ddots & \vdots \\
\dfrac{1}{Z_M-z_1} & \dfrac{1}{Z_M-z_2} & \dots & \dfrac{1}{Z_M-z_p}
    \end{bmatrix},\nonumber\\
lZ = \begin{bmatrix}
\log(Z_1-z_1) & \log(Z_1-z_2) & \dots & \log(Z_1-z_p) \\
\log(Z_2-z_1) & \log(Z_2-z_2) & \dots & \log(Z_2-z_p) \\
\vdots & \vdots & \ddots & \vdots \\
\log(Z_M-z_1) & \log(Z_M-z_2) & \dots & \log(Z_M-z_p)
    \end{bmatrix}.
\end{align}
The vector $b=(\tilde{u}(Z),\tilde{v}(Z))^T$ corresponds to the two boundary conditions on $\partial\Omega$. We compute the least-squares problem using the backslash command in MATLAB to obtain
\begin{equation}
x=[\mathrm{Re}\{a_j^f, b_j^f, c_{ij}^f, a_j^g, b_j^g, c_{ij}^g, d_i^f, d_i^g\},\mathrm{Im}\{a_j^f, b_j^f, c_{ij}^f, a_j^g, b_j^g, c_{ij}^g, d_i^f, d_i^g\}]^T,
\end{equation}
and thus all complex coefficients in $\hat{f}(z)$ and $\hat{g}(z)$, which satisfy the given boundary conditions in a least-squares manner.

After finding the Goursat functions, we construct function handles for physical quantities $u(z)$, $v(z)$, $p(z)$ and $\omega(z)$ using \cref{eq:bc1,eq:bc2}. The function handles enable the evaluation of any required physical quantity at any given point in $\Omega$.

\subsection{Application to other Stokes flow problems}
\label{sec:application_examples}

\begin{figure}[htbp]
  \centering
  \includegraphics[width=0.8\textwidth]{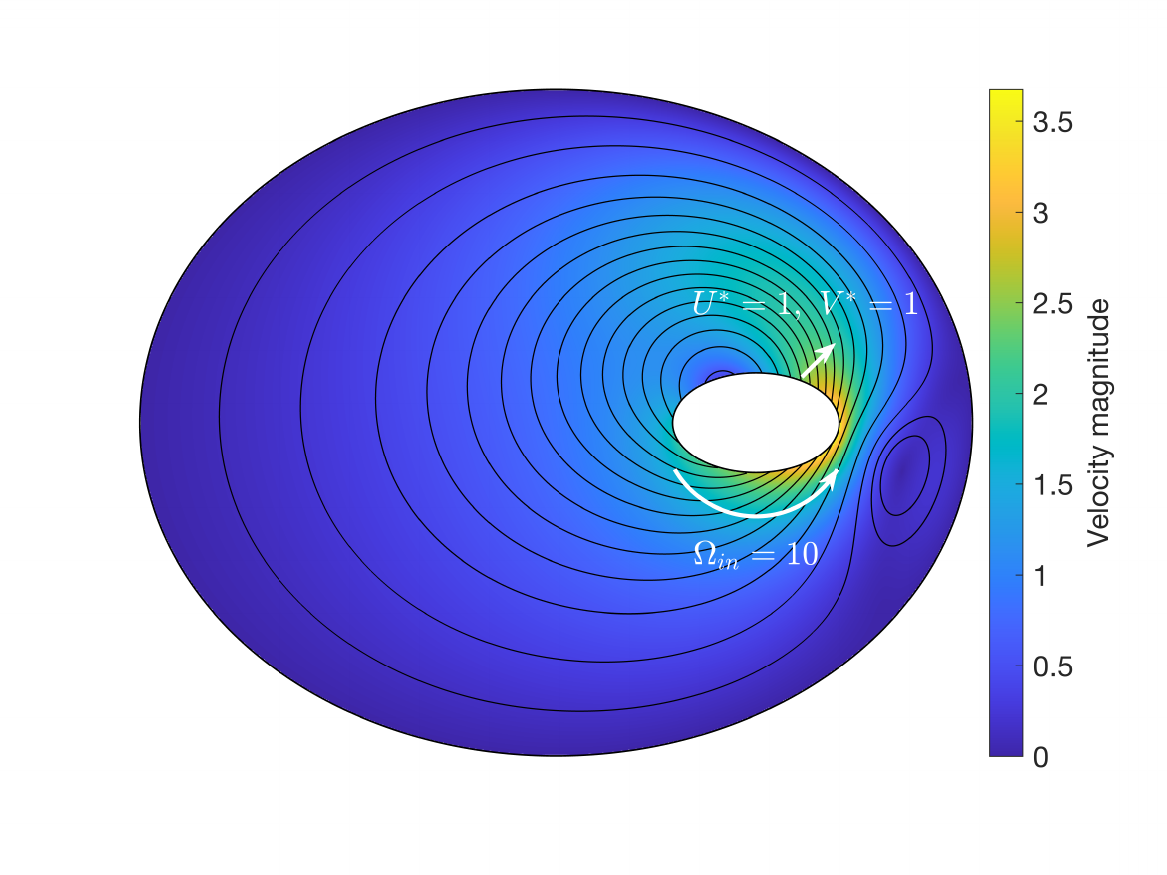}
  \caption{Stokes flow around a translating and rotating elliptical cylinder inside a fixed elliptical cylinder. The same parameter values as for Case `c' in \cref{tab:parameter_value} are used here. The outer ellipse has eccentricity $0.6$ and the inner ellipse has eccentricity $0.8$. The translation and rotation of the inner ellipse are indicated by white arrows.}
  \label{fig:multi_case1}
\end{figure}

\Cref{fig:multi_case1} shows 2D Stokes flow around a translating and rotating elliptical cylinder inside a fixed elliptical cylinder. This setting has potential biomedical applications to kidney stone removal problems \cite{Williams2020}. The parameter values for Case `c' in \cref{tab:parameter_value} are used here to set the boundary conditions, but with the outer cylinder replaced by an ellipse with eccentricity 0.6 and the inner cylinder by an ellipse with eccentricity 0.8, and $A_{out}$ and $A_{in}$ now representing the lengths of the semi-minor axes. The solution is computed to 7-digit accuracy with a degree 40 polynomial and a degree 120 Laurent series. Note that for non-circular holes, a higher degree Laurent series is usually required for good precision.

\begin{figure}[htbp]
  \centering
  \includegraphics[trim={0 4cm 0 4cm},clip,width=0.8\textwidth]{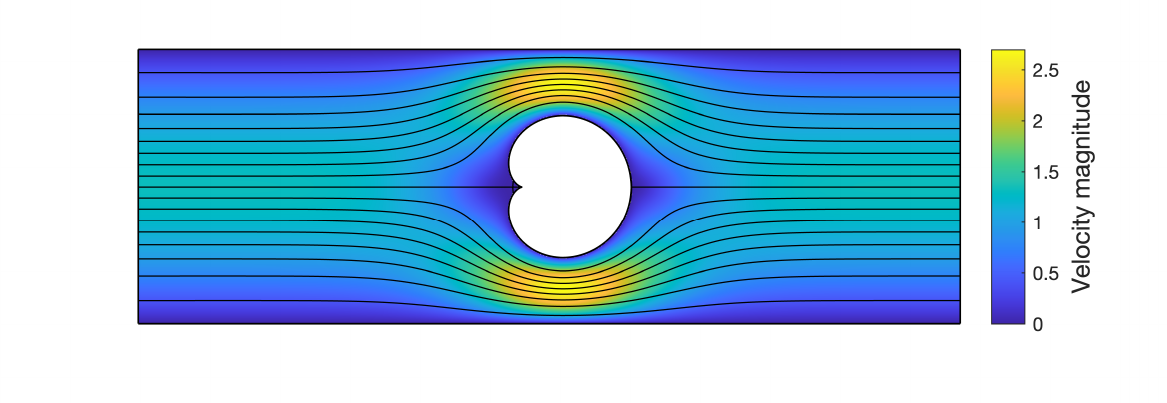}
  \caption{Stokes flow around a heart-shaped hole in a channel. The computation is carried out using a degree $120$ polynomial and a degree $80$ Laurent series.}
  \label{fig:multi_case2}
\end{figure}

\Cref{fig:multi_case2} shows the Stokes flow around a heart-shaped hole in a channel, illustrating how readily our solver can be extended to other shapes. Boundary conditions of constant pressure and parallel flow are imposed at the channel inlet and outlet, and a zero velocity condition is imposed on the walls. The solution is computed to 10-digit accuracy (12 to 13-digit accuracy along the hole boundary) with a degree 120 polynomial and a degree 80 Laurent series in 0.9 second. One can compute the solution to 6-digit accuracy with a degree 80 polynomial and a degree 40 Laurent series in 0.3 second.

\begin{figure}[htbp]
  \centering
  \includegraphics[width=0.55\textwidth]{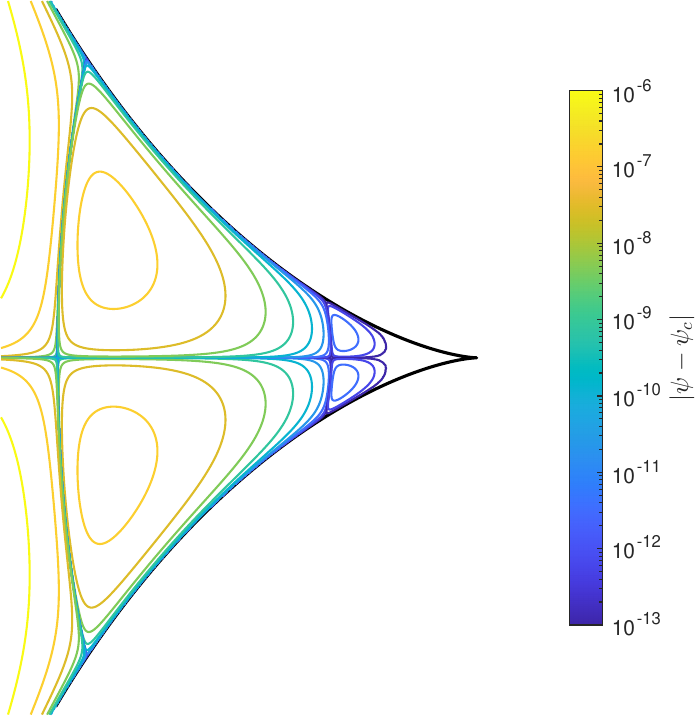}
  \caption{Moffatt eddies near the cusp of the heart-shaped hole in \cref{fig:multi_case2}. The colour scale represents the magnitude of the deviation of the stream function from $\psi_c$, the value along the centreline.}
  \label{fig:heart_closeup}
\end{figure}

\Cref{fig:heart_closeup} shows two pairs of Moffatt eddies \cite{Moffatt1964} near the cusp on the left side of the heart-shaped hole in \cref{fig:multi_case2}. If the stream function is set to 0 along the centreline, the first pair have stream functions on the order of $10^{-7}$ and the second on the order of $10^{-12}$. The third pair in this theoretically infinite series of eddies would be at a level near machine precision. This example demonstrates the great accuracy of our methods. For this problem, it takes 45 microseconds per point to compute a stream function or velocity.

\begin{figure}[htbp]
  \centering
  \includegraphics[trim={0 1cm 0 1cm},clip,width=\textwidth]{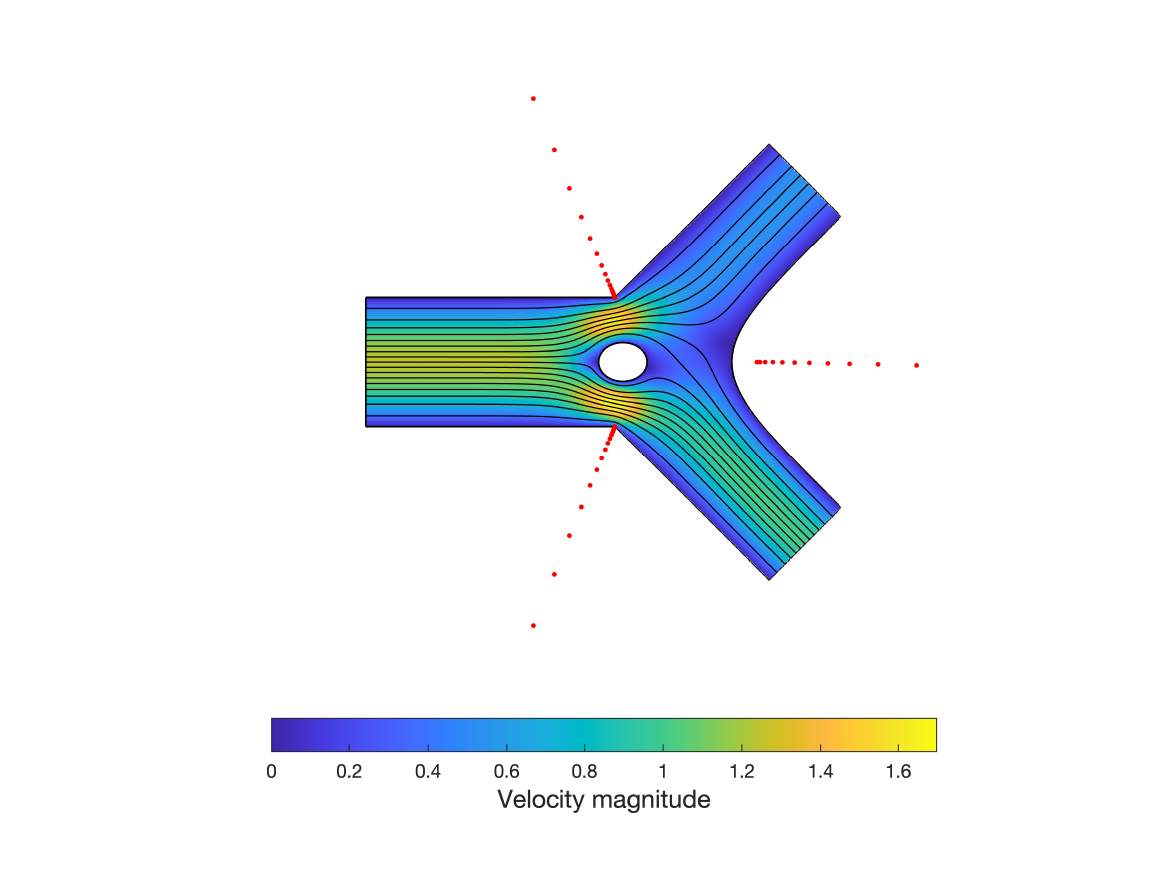}
  \caption{Stokes flow around a steady ellipse in a bifurcation with a smooth corner and two sharp corners. The computation is carried out to $6$-digit accuracy using a degree $96$ polynomial, a degree $48$ Laurent series, $48$ poles exponentially clustered near each sharp corner and AAA poles near the smooth corner.}
  \label{fig:bifurcation_particle}
\end{figure}

We conclude this section with an example case that combines all the methods we have introduced: the lightning solver for sharp corners \cite{Gopal2019a,Brubeck2022}, the AAA rational approximation for smooth boundaries \cite{Nakatsukasa2018,Costa2023} and the series method for multiply connected domains \cite{Axler1986,Finn2003,Trefethen2018}. This is an application of the complete LARS algorithm. \Cref{fig:bifurcation_particle} shows the Stokes flow around a steady ellipse (or an elliptical hole) within a bifurcation. The bifurcation has two sharp corners, where poles are exponentially clustered, and a smooth corner, where poles are placed using the AAA algorithm. Boundary conditions of constant pressure and parallel flow are imposed on the channel inlet and outlets, but the upper branch has a higher outlet pressure than that of the lower branch. Zero velocity conditions are imposed on both the channel walls and the ellipse boundary.

For sharp corners, poles are exponentially clustered along the exterior bisectors, as shown in previous work \cite{Brubeck2022}. For the smooth corner, AAA rational approximation places the poles outside the corner boundary using the boundary function $F=\bar{Z}$. It is interesting to note that the AAA algorithm clusters poles towards a branch point along the exterior bisector of the smooth boundary automatically in 0.03 second, without any knowledge of the singularities near that geometry. A degree 48 Laurent series was added for the elliptical hole (or a 2D steady elliptical particle) with corresponding logarithmic terms. Next, we carried out VA orthogonalization \cite{Brubeck2021} for a degree 96 polynomial, the Laurent series and three sets of poles near three corners. As the last step, we solved a linear least-squares problem using the backslash command in MATLAB to achieve a solution to 6-digit accuracy. The entire computation takes 1.5 seconds on a standard laptop. In engineering applications, e.g., for microparticle transport problems \cite{Audet1987,Doyeux2011} when Stokes flow needs to be simulated at multiple time steps, one can obtain a reasonably accurate solution (losing 1-digit accuracy) in a much shorter time by reducing the degrees of polynomial and series.

\section{Discussion}

In this paper, we have presented LARS, an algorithm that uses lightning and AAA rational approximation to compute 2D Stokes flows in bounded domains. The algorithm uses lightning approximation for sharp corners \cite{Gopal2019a,Brubeck2022}, AAA rational approximation for curved boundaries \cite{Nakatsukasa2018,Costa2023} and a series method for holes \cite{Axler1986,Trefethen2018}. Vandermonde with Arnoldi orthogonalization \cite{Brubeck2021} is carried out for each part to construct a well-conditioned basis, and the Goursat functions for 2D Stokes flows \cite{Goursat1898} are approximated by solving a linear least-squares problem.

One advantage of LARS is its great speed and accuracy, thanks to the root-exponential convergence of the lightning algorithm \cite{Gopal2019a}. We can compute a solution to a few digits of accuracy by just placing a few poles exponentially clustered near the singular corners or branch points near curved boundaries. The application of the lightning algorithm to Stokes flows in polygons has been presented in previous work \cite{Brubeck2021}. Similar speed and accuracy have been shown in this work when computing Stokes flows in curved boundary domains using AAA poles, which are also clustered towards the branch point. The beauty of using lightning and AAA rational approximation to compute Stokes flows is that it allows a rational function to capture the singularities using poles located outside the domain, while preserving analyticity inside the domain. The advantage of rational approximation over polynomial approximation in computing Stokes flows was shown in the constricted channel example in \cref{sec:aaa}.

The LARS solver is also suitable for quasi-steady computations with moving boundaries, where the boundary position and velocity are updated at each time step based on the flow field. This is because it is fast enough to perform hundreds of simulations in a few minutes. Since this algorithm does not require the domain discretisation at each time step, it saves both computer memory and computation time. We believe the proposed solver has many more potential applications to time-dependent problems than the steady-state problems presented in \cref{sec:solver_application}. For example, we are currently using it to investigate the transport of microparticles in channels with bifurcations \cite{Doyeux2011} with applications to the transport of microthrombi in the human cerebral microvasculature \cite{Xue2022}.

It should be noted that the algorithm presented is not limited to computing Stokes flows, although this is the focus of our applications. It can be applied to other biharmonic problems, e.g. the vibration of plates in solid mechanics \cite{Love1927}, with changes only in the boundary conditions.

One limitation of LARS is that it is only applicable to 2D geometries, because it is based on Goursat representations and rational functions. In fact, this is true for most applied complex variables techniques \cite{Ablowitz2003}. For 3D Stokes flow problems, one can use other numerical methods mentioned including finite element methods \cite{Logg2012} or boundary integral methods \cite{Pozrikidis1992}.

In addition, LARS currently only works for flows in bounded domains. As discussed previously \cite{Brubeck2022}, rational approximations should be able to treat Stokes flow problems in unbounded domains. The computation of unbounded Stokes flows requires careful consideration of the boundary conditions at infinity. For periodic boundary problems, e.g. Stokes flow through a 2D channel with periodic boundary geometry, using trigonometric polynomials might lead to better results than using the Runge terms in \cref{eq:rational}. For example, see \cite{Baddoo2021b} for an extended AAA algorithm using a trigonometric barycentric formula to approximate periodic functions.

Lastly, there has been recent progress in lightning and AAA rational approximations that may lead to a better 2D Stokes flow solver. Baddoo and Trefethen \cite{Baddoo2021a} developed a log-lightning method that has been shown to have faster convergence than the original lightning method \cite{Gopal2019a} when computing Green's function for a square. Following this, a log-lightning Stokes solver could be developed in future work. In addition, a new study considers AAA approximation on a continuum \cite{Driscoll2023}, rather than the vector of discrete points used in the original AAA approximation paper \cite{Nakatsukasa2018}. This can potentially lead to faster and more robust approximations.

To conclude, we have shown that it is now possible to compute 2D Stokes flows in very general domains using rational approximation. The ``LARS'' algorithm is simple and easy to implement for a variety of 2D Stokes flow problems. The computation usually takes less than a second to obtain a solution to at least 6-digit accuracy.

\appendix
\section{MATLAB codes}
\label{sec:code}
The MATLAB codes for computing the example cases in \cref{sec:aaa,sec:series_method,fig:heart_closeup,fig:bifurcation_particle} have been posted in a GitHub repository at \url{https://github.com/YidanXue/LARS}.

\section*{Acknowledgments}
We are grateful for very helpful discussions with Stephen Payne and Howard Stone. We also thank two anonymous referees for their comments and suggestions.

\bibliographystyle{siamplain}
\bibliography{references}

\begin{thebibliography}{10}

\bibitem{Ablowitz2003}
{\sc M.~J. Ablowitz and A.~S. Fokas}, {\em Complex variables: introduction and applications}, Cambridge University Press, 2003.

\bibitem{Audet1987}
{\sc D.~Audet and W.~Olbricht}, {\em The motion of model cells at capillary bifurcations}, Microvasc. Res., 33 (1987), pp.~377--396.

\bibitem{Axler1986}
{\sc S.~Axler}, {\em Harmonic functions from a complex analysis viewpoint}, Amer. Math. Monthly, 93 (1986), pp.~246--258.

\bibitem{Baddoo2021b}
{\sc P.~J. Baddoo}, {\em The {AAAtrig} algorithm for rational approximation of periodic functions}, SIAM J. Sci. Comput., 43 (2021), pp.~A3372--A3392.

\bibitem{Baddoo2021a}
{\sc P.~J. Baddoo and L.~N. Trefethen}, {\em Log-lightning computation of capacity and {Green}'s function}, Maple Transactions, 1 (2021).

\bibitem{Ballal1976}
{\sc B.~Ballal and R.~Rivlin}, {\em Flow of a {Newtonian} fluid between eccentric rotating cylinders: inertial effects}, Arch. Rational Mech. Anal., 62 (1976), pp.~237--294.

\bibitem{Banjai2008}
{\sc L.~Banjai}, {\em Revisiting the crowding phenomenon in {Schwarz}--{Christoffel} mapping}, SIAM J. Sci. Comput., 30 (2008), pp.~618--636.

\bibitem{Brubeck2021}
{\sc P.~D. Brubeck, Y.~Nakatsukasa, and L.~N. Trefethen}, {\em Vandermonde with {Arnoldi}}, SIAM Rev., 63 (2021), pp.~405--415.

\bibitem{Brubeck2022}
{\sc P.~D. Brubeck and L.~N. Trefethen}, {\em Lightning {Stokes} solver}, SIAM J. Sci. Comput., 44 (2022), pp.~A1205--A1226.

\bibitem{Costa2023}
{\sc S.~Costa and L.~N. Trefethen}, {\em {AAA}-least squares rational approximation and solution of {Laplace} problems}, in European Congr. Math., {Hujdurovi\'{c} et al.}, ed., EMS Press, 2023.

\bibitem{Davis1958}
{\sc P.~Davis and H.~Pollak}, {\em On the analytic continuation of mapping functions}, Trans. Amer. Math. Soc., 87 (1958), pp.~198--225.

\bibitem{Doyeux2011}
{\sc V.~Doyeux, T.~Podgorski, S.~Peponas, M.~Ismail, and G.~Coupier}, {\em Spheres in the vicinity of a bifurcation: elucidating the {Zweifach}--{Fung} effect}, J. Fluid Mech., 674 (2011), pp.~359--388.

\bibitem{Driscoll2023}
{\sc T.~Driscoll, Y.~Nakatsukasa, and L.~N. Trefethen}, {\em {AAA} rational approximation on a continuum}, SIAM J. Sci. Comput., to appear (2023).

\bibitem{Driscoll2014}
{\sc T.~A. Driscoll, N.~Hale, and L.~N. Trefethen}, {\em Chebfun {Guide}}, 2014, \url{www.chebfun.org}.

\bibitem{Dvinsky1987a}
{\sc A.~Dvinsky and A.~S. Popel}, {\em Motion of a rigid cylinder between parallel plates in {Stokes} flow: Part 1: Motion in a quiescent fluid and sedimentation}, Comput. Fluids, 15 (1987), pp.~391--404.

\bibitem{Dvinsky1987b}
{\sc A.~Dvinsky and A.~S. Popel}, {\em Motion of a rigid cylinder between parallel plates in {Stokes} flow: Part 2: {Poiseuille} and {Couette} flow}, Comput. Fluids, 15 (1987), pp.~405--419.

\bibitem{Finn2001}
{\sc M.~D. Finn and S.~M. Cox}, {\em Stokes flow in a mixer with changing geometry}, J. Eng. Math., 41 (2001), pp.~75--99.

\bibitem{Finn2003}
{\sc M.~D. Finn, S.~M. Cox, and H.~M. Byrne}, {\em Chaotic advection in a braided pipe mixer}, Phys. Fluids, 15 (2003), pp.~L77--L80.

\bibitem{Fokas2003}
{\sc A.~Fokas and A.~Kapaev}, {\em On a transform method for the {Laplace} equation in a polygon}, IMA J. Appl. Math., 68 (2003), pp.~355--408.

\bibitem{Frazer1926}
{\sc R.~A. Frazer}, {\em {III}. {On} the motion of circular cylinders in a viscous fluid}, Proc. R. Soc. Lond. A, 225 (1926), pp.~93--130.

\bibitem{Gopal2019b}
{\sc A.~Gopal and L.~N. Trefethen}, {\em Representation of conformal maps by rational functions}, Numer. Math., 142 (2019), pp.~359--382.

\bibitem{Gopal2019a}
{\sc A.~Gopal and L.~N. Trefethen}, {\em Solving {Laplace} problems with corner singularities via rational functions}, SIAM J. Numer. Anal., 57 (2019), pp.~2074--2094.

\bibitem{Goursat1898}
{\sc E.~Goursat}, {\em Sur l'\'{e}quation ${\Delta}{\Delta}u=0$}, Bull. Soc. Math. France, 26 (1898), pp.~236--237.

\bibitem{Jeffery1922}
{\sc G.~B. Jeffery}, {\em The rotation of two circular cylinders in a viscous fluid}, Proc. R. Soc. Lond. A, 101 (1922), pp.~169--174.

\bibitem{Jeong2001}
{\sc J.-T. Jeong}, {\em Slow viscous flow in a partitioned channel}, Phys. Fluids, 13 (2001), pp.~1577--1582.

\bibitem{Kelley2023}
{\sc D.~H. Kelley and J.~H. Thomas}, {\em Cerebrospinal fluid flow}, Annu. Rev. Fluid Mech., 55 (2023).

\bibitem{Kruger2017}
{\sc T.~Kr{\"u}ger, H.~Kusumaatmaja, A.~Kuzmin, O.~Shardt, G.~Silva, and E.~M. Viggen}, {\em The lattice {Boltzmann} method}, Springer, 2017.

\bibitem{Lauga2009}
{\sc E.~Lauga and T.~R. Powers}, {\em The hydrodynamics of swimming microorganisms}, Rep. Prog. Phys., 72 (2009), p.~096601.

\bibitem{Logg2012}
{\sc A.~Logg, K.-A. Mardal, and G.~Wells}, {\em Automated solution of differential equations by the finite element method: The {FEniCS} book}, vol.~84, Springer, 2012.

\bibitem{Love1927}
{\sc A.~E.~H. Love}, {\em A treatise on the mathematical theory of elasticity}, Cambridge University Press, Cambridge, 4th~ed., 1927.

\bibitem{Luca2018a}
{\sc E.~Luca and D.~G. Crowdy}, {\em A transform method for the biharmonic equation in multiply connected circular domains}, IMA J. Appl. Math., 83 (2018), pp.~942--976.

\bibitem{Luca2018b}
{\sc E.~Luca and S.~G. Llewellyn~Smith}, {\em Stokes flow through a two--dimensional channel with a linear expansion}, Q. J. Mech. Appl. Math., 71 (2018), pp.~441--462.

\bibitem{Manga1996}
{\sc M.~Manga}, {\em Dynamics of drops in branched tubes}, J. Fluid Mech., 315 (1996), pp.~105--117.

\bibitem{Menikoff1980}
{\sc R.~Menikoff and C.~Zemach}, {\em Methods for numerical conformal mapping}, J. Comput. Phys., 36 (1980), pp.~366--410.

\bibitem{Moffatt1964}
{\sc H.~K. Moffatt}, {\em Viscous and resistive eddies near a sharp corner}, J. Fluid Mech., 18 (1964), pp.~1--18.

\bibitem{Nakatsukasa2018}
{\sc Y.~Nakatsukasa, O.~S{\`e}te, and L.~N. Trefethen}, {\em The {AAA} algorithm for rational approximation}, SIAM J. Sci. Comput., 40 (2018), pp.~A1494--A1522.

\bibitem{Ockendon1995}
{\sc H.~Ockendon and J.~R. Ockendon}, {\em Viscous flow}, Cambridge University Press, 1995.

\bibitem{Pozrikidis1992}
{\sc C.~Pozrikidis}, {\em Boundary integral and singularity methods for linearized viscous flow}, Cambridge University Press, 1992.

\bibitem{Price2003}
{\sc T.~Price, T.~Mullin, and J.~Kobine}, {\em Numerical and experimental characterization of a family of two-roll-mill flows}, Proc. R. Soc. A, 459 (2003), pp.~117--135.

\bibitem{Ranger1961}
{\sc K.~Ranger}, {\em A problem on the slow motion of a viscous fluid between two fixed cylinders}, Q. J. Mech. Appl. Math., 14 (1961), pp.~411--421.

\bibitem{Ranger1980}
{\sc K.~Ranger}, {\em Eddies in two dimensional {Stokes} flow}, Int. J. Eng. Sci., 18 (1980), pp.~181--190.

\bibitem{Secomb2017}
{\sc T.~W. Secomb}, {\em Blood flow in the microcirculation}, Annu. Rev. Fluid Mech., 49 (2017), pp.~443--461.

\bibitem{Tavakol2017}
{\sc B.~Tavakol, G.~Froehlicher, D.~P. Holmes, and H.~A. Stone}, {\em Extended lubrication theory: improved estimates of flow in channels with variable geometry}, Proc. R. Soc. A, 473 (2017), p.~20170234.

\bibitem{Trefethen2018}
{\sc L.~N. Trefethen}, {\em Series solution of {Laplace} problems}, ANZIAM J., 60 (2018), pp.~1--26.

\bibitem{Trefethen2020}
{\sc L.~N. Trefethen}, {\em Numerical conformal mapping with rational functions}, Comput. Methods Funct. Theory, 20 (2020), pp.~369--387.

\bibitem{Walsh1929}
{\sc J.~Walsh}, {\em The approximation of harmonic functions by harmonic polynomials and by harmonic rational functions}, Bull. Amer. Math. Soc., 35 (1929), pp.~499--544.

\bibitem{Wiener1931}
{\sc N.~Wiener and E.~Hopf}, {\em \"{U}ber eine {Klasse} singul\"{a}rer {Integralgleichungen}}, Sem.–Ber. Preuss. Akad. Wiss., 31 (1931), pp.~696--706.

\bibitem{Williams2020}
{\sc J.~Williams, B.~Turney, D.~Moulton, and S.~Waters}, {\em Effects of geometry on resistance in elliptical pipe flows}, J. Fluid Mech., 891 (2020), p.~A4.

\bibitem{Xue2022}
{\sc Y.~Xue, T.~Georgakopoulou, A.-E. van~der Wijk, T.~I. J{\'o}zsa, E.~van Bavel, and S.~J. Payne}, {\em Quantification of hypoxic regions distant from occlusions in cerebral penetrating arteriole trees}, PLoS Comput. Biol., 18 (2022), p.~e1010166.

\bibitem{Yeo2011}
{\sc L.~Y. Yeo, H.-C. Chang, P.~P. Chan, and J.~R. Friend}, {\em Microfluidic devices for bioapplications}, Small, 7 (2011), pp.~12--48.

\end{thebibliography}
\end{document}